\documentclass[12pt,reqno]{amsart}   % Standard LaTeX

\usepackage[T1]{fontenc}
\usepackage{amsmath,amssymb,amsxtra,amscd}
\usepackage{amsthm,NewSegalNotation}
\usepackage[all]{xy}

\def\LinksPfeil{{\unitlength 1em\begin{picture}(0,.1)
\put(1.4,.1){\vector(-1,0){1.5}}
\end{picture}}}
\def\invlimname{{\unitlength.1em
\raisebox{-2.7\unitlength}{\begin{picture}(15.5,9.5)(0,0)
\put(0,2.7){$\operatorname{lim}$} \put(.05,-.1){\LinksPfeil}
\end{picture}}}}
\def\Invlim{\mathop{\invlimname}}
\newcommand {\invlim}[1] {\Invlim_{#1}}

\newtheorem{theorem}{Theorem}[section]
\newtheorem{proposition}[theorem]{Proposition}
\newtheorem{lemma}[theorem]{Lemma}
\newtheorem{definition}[theorem]{Definition}
\newtheorem{example}[theorem]{Example}
\newtheorem{corollary}[theorem]{Corollary}

\newtheorem{maintheorem}{Theorem}

\newtheorem*{question*}{Question}

\theoremstyle{definition}
\newtheorem{remark}[theorem]{Remark}

\begin{document}

\title{A Segal Conjecture for $p$-completed classifying spaces}
\author{K\'ari Ragnarsson}
\date{February 13, 2007}
\thanks{The author was supported by EPSRC grant GR/S94667/01 during most of this work}
 \address{Department of Mathematical Sciences, University of Aberdeen, Aberdeen AB24 3UE, United Kingdom}
\email{kari@maths.abdn.ac.uk}

\begin{abstract}
We formulate and prove a new variant of the Segal Conjecture
describing the group of homotopy classes of stable maps from the
$p$-completed classifying space of a finite group $G$ to the 
classifying space of a compact Lie group $K$ as the $p$-adic
completion of the Grothendieck group $A_p(G,K)$ of finite principal \linebreak
$(G,K)$-bundles whose isotropy groups are $p$-groups. 
Collecting the result for different primes $p$, we get a new and 
simple description of the group of homotopy classes of stable 
maps between (uncompleted) classifying spaces of groups. 
This description allows us to determine the kernel of the map from
the Grothendieck group $A(G,K)$ of finite principal $(G,K)$-bundles
to the group of homotopy classes of stable maps from $BG$ to $BK$.
\end{abstract}

 \maketitle
\section*{Introduction}
Inspired by Atiyah's description \cite{At} of the representable
complex periodic $K$-theory of a finite group $G$ as the
completion of the complex representation ring $R[G]$ at its
augmentation ideal,
  \[\hat{R}[G] \stackrel{\cong}{\longrightarrow} KU^{0}(BG),\]
Segal conjectured that the zeroth stable cohomotopy group of a
finite group $G$ could analogously be described as the completion
of the Burnside ring $A(G)$ of isomorphism classes of virtual
finite $G$-sets (\cite{Dieck}), with respect to its augmentation
ideal $I(G)$,
\[
  A(G)^{\wedge}_{I(G)}  \stackrel{\cong}{\longrightarrow} \pi^0_S(BG_+).
\]

After many partial results by various authors, this conjecture was
eventually settled in the affirmative by Carlsson in \cite{Car}.
In this account we mention only the three major contributions
which in the end combined to form a solution of the conjecture,
but the reader is encouraged to look up \cite{Car} for an
account of the history of the Segal conjecture.
Adams-Gunawardena-Miller proved the conjecture for elementary
abelian $p$-groups in \cite{AGM}. May-McClure simplified the
conjecture in \cite{MM}, where they proved that if the conjecture
holds for all finite $p$-groups, then it holds for all finite
groups. The solution of the conjecture was completed by Carlsson
in \cite{Car}, where he supplied an inductive proof showing that
if the Segal conjecture is true for elementary abelian $p$-groups
then it is true for all finite $p$-groups.

The Segal conjecture was extended by Lewis--May--McClure in
\cite{LMM} to describe the group of homotopy classes of stable
maps between classifying spaces of finite groups, and later by May--Snaith--Zelewski in \cite{MSZ}
to allow the classifying space of a compact Lie group in the target.
For a finite
group $G$ and a compact Lie group $K$, they consider the Burnside
module $A(G,K)$. This is the Grothendieck group completion of isomorphism 
classes of finite principal $(G,K)$-bundles, which are principal $K$-bundles 
with finitely many orbits in the category of $G$-spaces. The Burnside module $A(G,K)$ 
is an $A(G)$-module which as a $\Z$-module is free
with one basis element $[H,\varphi]$ for each conjugacy class of
pairs $(H,\varphi)$ consisting of a subgroup \mbox{$H \leq G$} and
a group homomorphism \mbox{$\varphi \colon H \to K$}.
Letting $\PtdStableMaps{BG}{BK}$ denote the group of homotopy
classes of stable maps \mbox{$\Stable{{BG}_+} \to
\Stable{{BK}_+}$}, there is a natural homomorphism
 \[\alpha \colon A(G,K) \longrightarrow \PtdStableMaps{BG}{BK},\]
which sends a basis element $[H,\varphi]$ to the stable map
$\Stable{B\varphi_+} \circ tr_H$, where $tr_H$ is the transfer map 
discussed in Section \ref{sub:Transfer}. In \cite{LMM} and \cite{MSZ} it is shown
that the Segal conjecture implies that this map is a
completion with respect to $I(G)$. The original Segal conjecture
is the special case when $K$ is the trivial group.

The $I(G)$-adic completion of $A(G,K)$ is in general difficult to
describe, but in the special case where $G$ is a finite $p$-group
$S$, May-McClure \cite{MM} gave a simple description. Letting
$\widetilde{A}(S,K)$ denote the quotient module obtained from
$A(S,K)$ by quotienting out basis elements of the form $[P,\triv]$
where $\triv$ is a constant homomorphism, they noticed that the
$I(S)$-adic topology on $\widetilde{A}(S,K)$ coincides with the
$p$-adic topology, and therefore the induced map
 \[\widetilde{\alpha} \colon \widetilde{A}(S,K) \longrightarrow \StableMaps{BS}{BK}\]
is a $p$-adic completion.

The main result of this paper can be regarded as an extension of
this result of May-McClure, where we drop the condition that $G$ be
a $p$-group but $p$-complete $BG$ instead. In other words we
describe $\StableMaps{\pComp{BG}}{BK}$ in similarly simple terms
as a $p$-completion of a certain submodule of $\widetilde{A}(G,K)$.
Let $A_p(G,K)$ be the submodule of $A(G,K)$ generated by those 
$(G,K)$-bundles whose isotropy groups are $p$-groups, and 
let $\widetilde{A}_p(G,K)$ denote the corresponding submodule of 
$A_p(G,K)$. We prove the following theorem, which
appears in the text as Theorem \ref{thm:NewSegal}.
\begin{maintheorem} \label{mthm:NewSegal}
For a finite group $G$ and a compact Lie group $K$, the homomorphism
\[
  \widetilde{A}_p(G,K) \xrightarrow{\widetilde{\alpha}_p} \StableMaps{BG}{BK} \xrightarrow{\iota_p^*} \StableMaps{\pComp{BG}}{BK}
\]
induced by $\alpha$ is a $p$-adic completion.
\end{maintheorem}
Here \mbox{$\iota_p \colon \pComp{\ClSpectrum{G}} \hookrightarrow
\ClSpectrum{G}$} is the natural wedge summand inclusion obtained
from the natural splitting \mbox{$ \Stable{BG} \simeq \bigvee_q
\Stable{\qComp{BG}}$}, as discussed in Section \ref{sec:pComps}.

Since $\alpha$ is natural in $G$ and finite 
$K$ we can state the finite group version of this theorem in terms
of isomorphisms of categories. This is done as Corollary
\ref{Cor:NewSegalCat} after we have developed the necessary
framework.

As a consequence of the stable splitting of the classifying space of
a finite group into its $q$-completed components, where $q$ runs over
all primes, one easily obtains a new description of $\StableMaps{BG}{BK}$. 
The result is similar in spirit to Minami's description of $A(G)^{\wedge}_I$ 
for a compact Lie group $G$ in \cite{Min}. Indeed, when $K$ is the trivial 
group one recovers Minami's result for finite $G$.
This result is Theorem 
\ref{thm:NewSegalSum} in the text.
\begin{maintheorem} \label{mthm:NewSegalSum}
For a finite group $G$ and a compact Lie group $K$, the
homomorphism $\alpha$ induces an isomorphism of $\Z$-modules
 \[\bigoplus_p \pComp{\widetilde{A}_p(G,K)} \stackrel{\cong}{\longrightarrow} \StableMaps{BG}{BK}.\]
\end{maintheorem}
This description is much simpler than the $I(G)$-adic completion,
but it has the drawback that it is not obvious how to decompose an
element in $\StableMaps{BG}{BK}$ into $p$-completed parts in
$\pComp{\widetilde{A}_p(G,K)}$. This matter is taken up in Section
\ref{sec:Decomposition} where we describe the map
 \[\widetilde{A}(G,K) \longrightarrow \StableMaps{BG}{BK} \longrightarrow \StableMaps{\pComp{BG}}{\pComp{BK}} \stackrel{\cong}{\longrightarrow} \pComp{\widetilde{A}_p(G,K)}\]
for every prime $p$. Loosely speaking, this map sends an element
$X$ of $A(G,K)$ to an element of $\pComp{\widetilde{A}_p(G,K)}$
whose $H$-fixed point sets for subgroups \mbox{$H \leq K\times G$}
agree with the $H$-fixed point sets of $X$ if $H$ is a finite $p$-group,
and are empty otherwise. (This description should not be taken
too literally as elements of $A(G,K)$ and
$\pComp{\widetilde{A}_p(G,K)}$ are not actual $(K \times G)$-bundles
and do not have fixed point sets. A precise statement can be found in
Section \ref{sec:Decomposition}.) As a consequence we deduce the
following description of the kernel of $\alpha$. 
\begin{maintheorem} \label{mthm:Kernel}
Let $G$ be a finite group and $K$ a compact Lie group.
\begin{itemize}
  \item[(a)] The kernel of the map
\[ \alpha \colon A(G,K) \longrightarrow \StableMaps{\pComp{BG}}{BK} \]
consists of those virtual $(G,K)$-bundles $[X - Y]$ whose
fixed-point sets under the action of any finite $p$-subgroup
\mbox{$P \leq K\times G$} satisfy
\[ |W(P)\backslash X^P| = |W(P) \backslash Y^P|,\]
where $W(P) = N_{K\times G}(P)/P.$ 
  \item[(b)]
 The kernel of the map
\[ \alpha \colon A(G,K) \longrightarrow \StableMaps{BG}{BK} \]
consists of those virtual $(G,K)$-bundles $[X - Y]$ such that for every prime
$p$ and every $p$-subgroup \mbox{$P \leq K\times G$}, the
$P$-fixed-point sets satisfy
\[ |W(P)\backslash X^P| = |W(P) \backslash Y^P|.\]
\end{itemize}
\end{maintheorem}
A more informative version of this theorem appears in the text as 
Corollary \ref{cor:Kernel} after some relevant concepts and notation 
have been introduced. In particular, not all $p$-subgroups $P$ 
need be considered, and there is a justification why the quotients 
of the fixed-point sets appearing are finite sets. Note also that when $K$ 
is finite one does not need to take the quotient by the $W(P)$-action (cf. Remark \ref{rem:DropW}).
\smallskip

The Segal conjecture has been presented
here in its weak form. There is a stronger form of the conjecture,
also due to Segal, describing higher stable cohomotopy groups of
classifying spaces of finite groups as a completion of certain
equivariant stable cohomotopy groups, constructed by Segal in
\cite{Se}. The proofs of the Segal conjecture mentioned above are
in fact proofs of the stronger version. Indeed, Carlsson's
inductive argument would not have been possible without the
presence of higher homotopy groups. The Lewis-May-McClure
generalization of the conjecture also extends the stronger version
of the Segal conjecture to describe spectra of stable maps between
classifying spaces. We refer the reader to
\cite{LMM} for this description. This raises an obvious,
interesting question. Namely, whether the results in this paper
can be extended to also describe higher homotopy groups of spectra
of stable maps between $p$-completed classifying spaces.

\smallskip
This paper is divided into five sections. In Section
\ref{sec:Preliminaries} we discuss some background material
necessary for the main discussion. In Section
\ref{sec:Subconjugacy} we introduce subconjugacy, which gives a
convenient filtration of Burnside modules. In Section
\ref{sec:NewSegal} we state and prove the new variants of the
Segal conjecture for $p$-completed classifying spaces. Alternative formulations of the main results, arising from different ways to overcome a technical nuisance involving basepoints, are discussed in Section \ref{sec:Alternative}. Finally, in Section \ref{sec:Decomposition} we
discuss how to decompose a stable map into its $p$-completed
components and describe the kernel of the $I(G)$-adic completion
map.

\section{Preliminaries} \label{sec:Preliminaries}
In this section we discuss background material which will be
needed in later sections. We begin by giving an overview of
notation and conventions which will be used throughout this paper.

Unless otherwise specified, $p$ is a fixed prime and all
cohomology is taken with $\Fp$-coefficients. For a space or
spectrum $X$, we let $\pComp{X}$ denote the Bousfield-Kan
$p$-completion \cite{Monster}, and for a $\Z$-module $M$, we let
$\pComp{M}$ denote the $p$-adic completion of $M$. Since we only
consider finitely generated modules in this paper, 
we have \mbox{$M \cong M \otimes \Zp$},
where $\Zp$ denotes the $p$-adic integers.

For a space $X$, let $X_+$ denote the pointed space obtained by
adding a disjoint basepoint to $X$. We use the shorthand notation
\mbox{$\PtdStable{X} := \Stable{X_+}$} and recall that there is a natural 
equivalence
\mbox{$\PtdStable{X} \simeq \Stable{X} \vee \SphereSpectrum$},
where $\SphereSpectrum$ denotes the sphere spectrum. All stable
homotopy will take place in the homotopy category of spectra,
which we denote by $\SpectraCat$.

Let $\Gr$ denote the category of finite groups. We will use
$\iota_H$ to denote a subgroup inclusion \mbox{$H \hookrightarrow
G$}, specifying the supergroup $G$ when there is danger of
confusion. Conjugations will come up frequently. For a group element
\mbox{$g \in G$} we let $c_g$ denote the conjugation isomorphism
\mbox{$x \mapsto gxg^{-1}$}. For a subgroup \mbox{$H \leq G $}, we
let $\lsup{H}{g}$ denote $c_g(H)$, and $H^g$ denote $c_g^{-1}(H)$.

We use the shorthand notations \mbox{$\ClSpectrum :=
\Stable{B(-)}$} and \mbox{$\PtdClSpectrum := \PtdStable(B(-))$},
regarded as functors \mbox{$\Gr \to \SpectraCat$}. Since
\mbox{$\Stable{(\pComp{BG})} \simeq \pComp{(\Stable{BG})}$} for a
finite group $G$, we will write $\pComp{\ClSpectrum{G}}$ without
danger of confusion. We denote by $\StableGr$ the category whose
objects are the finite groups and whose morphisms are homotopy
classes of stable maps between classifying spaces,
 \[\Mor{G_1}{G_2}{\StableGr} = \StableMaps{BG_1}{BG_2}.\]
Similarly we let $\CompStableGr$ be the category with the same
objects, but whose morphisms are homotopy classes of stable maps
between $p$-completed classifying spaces.
\[
  \Mor{G_1}{G_2}{\CompStableGr} = \CompStableMaps{(BG_1)}{(BG_2)}.
\]

\subsection{Bousfield-Kan $p$-completion of $\ClSpectrum{G}$}\label{sec:pComps}
In this section we list the basic properties of Bousfield-Kan
completion of spectra at a prime. We recall how the suspension
spectrum of the classifying space of a finite group decomposes as a
wedge sum of its $p$-completions as $p$ runs over all primes and see
that this splitting is natural.

Bousfield-Kan completion at a prime $p$ is an endofunctor
\mbox{$\pComp{(-)}$} defined on either the category of spaces or
spectra, depending on the context. In either case its defining
property is that for a map \mbox{$f\colon X \to Y$}, the
$p$-completed map \mbox{$\pComp{f} \colon \pComp{X}\to\pComp{Y}$} is
a weak homotopy equivalence if and only if $f$ induces an
isomorphism in homology with $\Fp$ coefficients.

The $p$-completion functor comes with a natural transformation
\mbox{$\eta_p \colon Id \Rightarrow \pComp{(-)}$}. We say that a
space or spectrum $X$ is \emph{$p$-complete} if the map
\mbox{$\eta_p \colon X \to \pComp{X}$} is a weak equivalence. We
say that $X$ is \emph{$p$-good} if $\pComp{X}$ is $p$-complete.
Classifying spaces of finite groups and their suspension spectra
are $p$-good. Classifying spaces of finite $p$-groups and their
suspension spectra are $p$-complete.

As noted earlier, we have \mbox{$\pComp{(\Stable{BG})} \simeq
\Stable{(\pComp{BG})} $} for a finite group $G$, and so we can
denote both of these by $\pComp{\ClSpectrum{G}}$ without danger of
confusion. However, $p$-completion does not commute with
suspension in general, as we have
\[ \Stable{(\pComp{(BG_+)})} \simeq \Stable{(\pComp{BG})_+} \simeq \SphereSpectrum \vee \Stable{\pComp{BG}} \]
while
\[ \pComp{(\Stable{BG_+})} \simeq \pComp{\SphereSpectrum} \vee \pComp{\Stable{BG}}. \]
This difference is one of the underlying reasons for the basepoint
issues one has to circumvent in the formulation of the
$p$-completed Segal conjectures.

When $G$ is a finite group, $\ClSpectrum{G}$ is a torsion spectrum,
so by Sullivan's arithmetic square the natural maps \mbox{$\eta_q
\colon \ClSpectrum{G} \to \qComp{\ClSpectrum{G}}$} induce a natural
homotopy equivalence
\[
 h := \bigvee_q \eta_q \colon \ClSpectrum{G} \stackrel{\simeq}{\longrightarrow} \bigvee_q
 \qComp{\ClSpectrum{G}},
\]
where the wedge sum runs over all primes $q$. We obtain a natural
inclusion of $\pComp{\ClSpectrum{G}}$ as a wedge summand of
$\ClSpectrum{G}$, defined as the composite
\[
 \iota_p \colon \pComp{\ClSpectrum{G}} \hookrightarrow \bigvee_q \qComp{\ClSpectrum{G}} \xrightarrow{h^{-1}} \ClSpectrum{G}.
\]
with left homotopy inverse $\eta_p$.

\begin{lemma} For a finite group $G$ and a compact Lie group $K$, the 
$p$-completion functor sends $f \in \StableMaps{BG}{BK}$ to $\eta_p\circ f \circ \iota_p \in \CompStableMaps{BG}{BK}$.
\end{lemma}
\begin{proof} By naturality of $\eta_p$ we have $\pComp{f}\circ\eta_p = \eta_p\circ f$,
from which it follows that $\pComp{f} \simeq \pComp{f}\circ\eta_p \circ \iota_p= \eta_p\circ f \circ \iota_p.$ 
\end{proof}

Bousfield--Kan $p$-completion of spectra can also be thought of as Bousfield localization with respect to the homology theory $H_*(-,\Fp)$. Bousfield localization is discussed in \cite{Bou}. Taking this point of view, one of Bousfield's results in \cite{Bou} is that when $X$ is a connective spectrum, the $p$-completion of $X$ is homotopy equivalent to the function spectrum $F(S^{-1}\Z/p^{\infty},X)$, where $S^{-1}\Z/p^{\infty}$ is the desuspension of the Moore spectrum $S\Z/p^{\infty}$, with $Z/p^{\infty} = \Z[1/p]/\Z$.

\begin{lemma} Let $X$ and $Y$ be spectra such that $Y$ and $F(X,Y)$ are both connective. Then there is a natural homotopy equivalence
\[ \Ad \colon F(X,\pComp{Y}) \stackrel{\simeq}{\longrightarrow} \pComp{F(X,Y)} \]
\end{lemma}
\begin{proof} The map is just the adjunction
 \[F(X,F(S^{-1}\Z/p^{\infty},Y)) \cong F(S^{-1}\Z/p^{\infty},F(X,Y)). \]
\end{proof}

In this paper we study stable maps from the classifying space of a finite group to the classifying space of a compact Lie group.
These exhibit very unusual behaviour under $p$-completion as seen in the following lemma. 
%In the statement $F(X,Y)$ denotes the function
%spectrum of two spectra $X$ and $Y$, whose zeroth homotopy group is $[X,Y]$.
\begin{lemma} \label{lem:pCompsAgree}
For a finite group $G$ and a compact Lie group $K$ the maps in the sequence below are all homotopy equivalences.
\[ F(\pComp{\ClSpectrum{G}},\ClSpectrum{K}) \xrightarrow{\eta_p \circ -}  F(\pComp{\ClSpectrum{G}},\pComp{\ClSpectrum{K}}) \xrightarrow{- \circ \eta_p} F(\ClSpectrum{G},\pComp{\ClSpectrum{K}}) \xrightarrow{\Ad} \pComp{F(\ClSpectrum{G},\ClSpectrum{K})} \]
In particular they induce isomorphisms of groups of homotopy classes of maps.
\end{lemma}
\begin{proof} 
We obtain the homotopy equivalence $\Ad$ because $\ClSpectrum{K}$ is connective by construction as a suspension spectrum, and $F(\ClSpectrum{G},\ClSpectrum{K})$ is connective by \cite{LMM}. The map $- \circ \eta_p$ is a homotopy equivalence since the cofibre of $\eta_p$ is $H_*(-,\Fp)$-acyclic (\cite{Bou}).

The map $\eta_p \circ -$ factors as the composition
\[ F(\pComp{\ClSpectrum{G}},\ClSpectrum{K}) \xrightarrow{\eta_p} \pComp{F(\pComp{\ClSpectrum{G}},\ClSpectrum{K})} \xrightarrow{\Ad}   F(\pComp{\ClSpectrum{G}},\pComp{\ClSpectrum{K}}) \] 
Let $S$ be a Sylow subgroup of $G$. As explained in the next subsection, 
the spectrum $F(\pComp{\ClSpectrum{G}},\ClSpectrum{K})$ is a retract of the spectrum 
$F(\ClSpectrum{S},\ClSpectrum{K})$, which is $p$-complete and connective by \cite{LMM}. Hence $F(\pComp{\ClSpectrum{G}},\ClSpectrum{K})$ is \mbox{$p$-complete} and connective, and so the two maps in the factorization are both homotopy equivalences.
\end{proof}

When working with ``pointed'' classifying spectra one similarly has a sequence of homotopy equivalences
 \[  F(\pComp{(\PtdClSpectrum{G})},\pComp{(\PtdClSpectrum{K})}) \xrightarrow{- \circ \eta_p} F(\PtdClSpectrum{G},\pComp{(\PtdClSpectrum{K})}) \xrightarrow{\Ad} \pComp{F(\PtdClSpectrum{G},\PtdClSpectrum{K})}. \]
The function spectrum $F(\pComp{(\PtdClSpectrum{G})},\PtdClSpectrum{K})$ is no longer homotopy equivalent because it contains an uncompleted wedge summand $F(\pComp{\SphereSpectrum},\SphereSpectrum)$ (but $p$-completing this summand makes it homotopy equivalent).

\subsection{Transfers} \label{sub:Transfer}
Recall (see for example \cite{Ad}) that a subgroup inclusion
\mbox{$H \hookrightarrow G$}, has a transfer map
 \[tr_H \colon \PtdClSpectrum{G} \longrightarrow \PtdClSpectrum{H}\]
such that the composition $\PtdClSpectrum{\iota_H} \circ tr_H$
acts as multiplication by $[G:H]$ in singular cohomology (with any
coefficients). We will for the most part prefer to work with the
reduced transfer map \mbox{$\ClSpectrum{G} \to \ClSpectrum{H}$},
which we also denote by $tr_H$ to reduce notation. The composition
$\ClSpectrum{\iota} \circ tr_H$ involving the reduced transfer
also acts as multiplication by $[G:H]$ in singular cohomology.

When $S$ is a Sylow subgroup of $G$, the index $[G:S]$ is a unit
in $\Fp$ and therefore \mbox{$\ClSpectrum{\iota_S} \circ tr_S$}
induces an isomorphism on $\Coh{\ClSpectrum{G}} = \Coh{BG}$. After
$p$-completion, $\ClSpectrum{\iota_S} \circ tr_S$ consequently
becomes a homotopy equivalence
 \[\pComp{\ClSpectrum{G}}\stackrel{\pComp{tr}}{\longrightarrow} \ClSpectrum{S} \stackrel{\pComp{\ClSpectrum{\iota}}}{\longrightarrow}  \pComp{\ClSpectrum{G}}\]
and so the $p$-completions of $\ClSpectrum{\iota_S}$ and $tr_S$
make $\pComp{\ClSpectrum{G}}$ a wedge summand of
\mbox{$\pComp{\ClSpectrum{S}} \simeq \ClSpectrum{S}$}.

When $H$ is a subgroup of $G$, we will often denote the
$p$-completed transfer $\pComp{(tr_H)}$ by $tr_H$. Similarly, when
\mbox{$\varphi \colon H \to G$} is a group homomorphism,
we will often denote the $p$-completed map
$\pComp{(\ClSpectrum{\varphi})}$ by $\ClSpectrum{\varphi}$. This
is done to reduce notation. As we always distinguish clearly
between $\pComp{\ClSpectrum{G}}$ and $\ClSpectrum{G}$ (except when
$S$ is a $p$-group and the spectra are homotopy equivalent
anyway), this should not cause any confusion.

\subsection{Linear Categories}
For a ring $R$, an \emph{$R$-linear category} is a category in which
the morphism sets are $R$-modules and composition is bilinear. It is
perhaps more common to refer to $\Z$-linear categories as
\emph{pre-additive} categories, but we use different terminology to
emphasize the role of the ground ring $R$. We say that a functor
\mbox{$F \colon \mathcal{C} \to \mathcal{D}$} is a
\emph{functor of $R$-linear categories} if the categories
$\mathcal{C}$ and $\mathcal{D}$ are $R$-linear and for every pair of
objects \mbox{$c,c' \in \mathcal{C}$} the map
 \[F \colon \Mor{c}{c'}{\mathcal{C}} \to \Mor{F(c)}{F(c')}{\mathcal{D}} \]
is a morphism of $R$-modules. An \emph{isomorphism of $R$-linear
categories} is a functor of $R$-linear categories that is also an
isomorphism of categories.

\begin{example}
If $R$ is a ring, then we can regard $R$ as an $R$-linear category
with one object $\circ_R$ such that
 \[\Mor{\circ_R}{\circ_R}{} = R.\]
We will sometimes take this point of view for the rings $\Z$ and
$\Zp$.
\end{example}

\subsection{The Burnside category}
For a finite group $G$ and a compact Lie group $K$, let
$\Morita(G,K)$ be the set of isomorphism classes of finite
principal $(G,K)$-bundles over finite $G$-sets. 
A $(G,K)$-bundle here means a $K$-bundle in the category of 
$G$-spaces, finiteness means the base space is finite, and principality means
that the $K$-action is free.
When $K$ is finite these are just $K$-free, finite $(K \times G)$-sets. 
The set $\Morita(G,K)$ is an abelian monoid under the
disjoint union operation and so we can consider its Grothendieck
group completion which we call the \emph{Burnside module of $G$
and $K$}, and denote by $A(G,K)$. The Burnside modules are abelian
groups, and hence $\Z$-modules, by construction. Their module
structure is well understood and will be described below after a
preliminary definition.

\begin{definition}
Let $G$ be a finite group and $K$ a compact Lie group. A
\emph{$(G,K)$-pair} is a pair $(H,\varphi)$ consisting of a subgroup
\mbox{$H \in G$} and a homomorphism \mbox{$\varphi \colon H
\to K$}. We say two $(G,K)$-pairs $(H,\varphi)$ and $(H',\varphi')$
are \emph{$(G,K)$-conjugate} if there exist elements \mbox{$g \in
G$} and \mbox{$k \in K$} such that \mbox{$c_g(H) = H'$} and the
following diagram commutes:
\[
\begin{CD}
{H} @ > {\varphi} >> {K} \\
@ V \cong V c_g V @ VV c_k V \\
{H'} @> {\varphi'} >> {K.} \\
\end{CD}
\]
Let $C(G,K)$ denote the set of conjugacy classes of $(G,K)$-pairs.
\end{definition}
When there is no danger of confusion we say
\emph{conjugacy} instead of $(G,K)$-conjugacy. We denote the
conjugacy class of a $(G,K)$-pair $(H,\varphi)$ by
$[H,\varphi]_{G}^{K}$, or just $[H,\varphi]$ when there is no danger
of confusion. 

\begin{lemma} \label{lem:pCompFactors}
For a finite group $G$ and a compact Lie group $K$, the
number of conjugacy classes of $(G,K)$-pairs is finite.
\end{lemma}
\begin{proof} One proves this easily using the folklore 
result (see for example \cite[Theorem 2.3 (1)]{MaPr:CptStCl} for a proof)
that for any finite group $F$ there are only finitely many 
conjugacy classes of subgroups of $K$ isomorphic to $F$. 
\end{proof}

A $(G,K)$-pair $(H,\varphi)$ gives rise to a finite principal $(G,K)$-bundle  
 \[(K\times G)/\Graph{H}{\varphi} \to G/H, \]
where
  \begin{equation} \label{eq:GraphDef}
    \Graph{H}{\varphi}  := \{ (\varphi(h),h) \mid h \in H \} \leq
    K\times G
  \end{equation}
is the (transposed) graph of $\varphi \colon H \to K$.  Two
$(G,K)$-pairs give rise to isomorphic $(G,K)$-bundles if and only if they 
are conjugate.

The ring $U(K\times G)$ of finite $(K\times G)$-complexes is introduced
in \cite{Dieck} where a basis over $\Z$ is also described. One can identify 
$A(G,K)$ with the submodule of $U(G \times K)$ generated by those 
$(K \times G)$-complexes whose isotropy groups are of the form $\Graph{H}{\varphi}$. 
This leads to the following description of $A(G,K)$.

\begin{proposition}  \cite{Dieck} \label{prop:BurnsideStructure}
For a finite group $G$ and a compact Lie group $K$, the Burnside
module $A(G,K)$ is a finitely generated, free $\Z$-module with one
basis element for each conjugacy class of $(G,K)$-pairs.
\end{proposition}
\begin{proof} The ring $U(K \times G)$ is free as a $\Z$-module with one basis
element for every conjugacy class of subgroups of $K \times G$. The Burnside module
$A(G,K)$ can be identified with the submodule generated by the subgroups of 
the form $\Graph{H}{\varphi}$ where $(H,\varphi)$ is a $(G,K)$-pair.
\end{proof}
By a slight abuse of notation we also denote the basis element corresponding to 
a $(G,K)$-pair $(H,\varphi)$ by $[H,\varphi]_{G}^{K}$ (or $[H,\varphi]$). 
For each basis element $[H,\varphi]$ we define a morphism
 \[\Coeff_{[H,\varphi]} \colon A(G,K) \longrightarrow \Z \]
by letting
$\Coeff_{[H,\varphi]}(X)$ be the coefficient at $[H,\varphi]$
in the basis decomposition of $X$. In other words we demand
that
 \[X = \sum_{[H,\varphi]} \Coeff_{[H,\varphi]}(X)\cdot [H,\varphi]\]
for all \mbox{$X \in A(G,K)$}. When appropriate, we will
also let $\Coeff_{[H,\varphi]}$ denote the analogous homomorphism
\mbox{$\pComp{A(G,K)} \to \Zp$}.

The Burnside modules
form the morphism sets of a certain linear category as described below.
For finite groups $G_1$ and $G_2$, and a compact Lie group $K$, there is a pairing
 \[\Morita(G_2,K) \times \Morita(G_2,G_1) \longrightarrow \Morita(G_1,K), \]
given by
 \[(X,Y) \longmapsto G_2 \backslash (X \times Y),\]
where $G_2$ acts by the diagonal. This extends to a bilinear pairing
 \[A(G_2,K) \times A(G_1,G_2) \longrightarrow A(G_1,K),\]
suggestively denoted by
 \[(X,Y) \longmapsto X \circ Y.\]
This pairing can be described on basis elements by the double
coset formula
 \[ [H_2,\varphi_2] \circ [H_1,\varphi_1] = \DCF{H_2}{\varphi_2}{H_1}{\varphi_1}{G_2}.\]

It is easy to see that this composition pairing satisfies the
associativity law and therefore we can make the following definition.
\begin{definition} 
The \emph{Burnside category} is the
$\Z$-linear category $\Burnside$ whose objects are the finite
groups and whose morphism sets are the Burnside modules,
 \[\Mor{G_1}{G_2}{\Burnside} = A(G_1,G_2),\]
where composition is given by the pairing described above.

Similarly, the  \emph{$p$-completed Burnside category} is the
$\Zp$-linear category $\CompBurnside$ whose objects are the finite
groups and whose morphism sets are the $p$-completed Burnside
modules,
 \[\Mor{G_1}{G_2}{\CompBurnside} = \pComp{A(G_1,G_2)}.\]
\end{definition}

Since we are focusing on $p$-local properties it will suit us, for
a finite group $G$ and a compact Lie group $K$, to study the submodule $A_p(G,K)$ of
the Burnside module $A(G,K)$ obtained by considering only
$(G,K)$-bundles whose isotropy groups are all $p$-groups.
Alternatively, this is the submodule generated by those basis
elements $[P,\varphi]$, where $P$ is a $p$-group. By the double coset formula, 
we see that these modules are
preserved by the composition pairing as described in the following
lemma.
\begin{lemma}\label{lem:pBurnsidePreserved}
For finite groups $G_1$ and $G_2$, and a compact Lie group $K$, we
have
  \[A(G_2,K) \circ A_p(G_1,G_2) \subset A_p(G_1,K).\]
In particular
  \[A_p(G_2,K) \circ A_p(G_1,G_2) \subset A_p(G_1,K).\]
\end{lemma}

Therefore we can make the following definition.

\begin{definition}
The \emph{$p$-isotropy Burnside category} is the $\Z$-linear
category $\Burnsidep$ whose objects are the finite groups and
whose morphism sets are the $p$-Burnside modules,
 \[\Mor{G_1}{G_2}{\Burnsidep} = A_p(G_1,G_2).\]

Similarly, the  \emph{$p$-completed $p$-Burnside category} is the
$\Zp$-linear category $\CompBurnsidep$ whose objects are the
finite groups and whose morphism sets are the $p$-completed
$p$-Burnside modules,
 \[\Mor{G_1}{G_2}{\CompBurnsidep} = \pComp{A_p(G_1,G_2)}.\]
\end{definition}

We conclude this discussion by introducing an augmentation for the Burnside
category. This will prove very useful for understanding the filtration of Burnside
modules introduced in Section \ref{sec:Subconjugacy}.
\begin{definition}
For a finite group $G$ and a compact Lie group $K$, the
\emph{orbit augmentation} of $A(G,K)$ is the homomorphism
\[ \Orbits \colon A(G,K) \longrightarrow \Z \]
obtained as the group completion of the monoid morphism
\[ \Morita(G,K) \longrightarrow \Z,~ X \longmapsto |K \backslash X|.\]
The \emph{orbit augmentation} of $\pComp{A(G,K)}$ is the
homomorphism
\[ \Orbits \colon \pComp{A(G,K)} \to \Zp\]
obtained upon $p$-completion.
\end{definition}
When $K$ is finite we have \mbox{$ \Orbits(X) =
|X|/|K|,$} since $X$ is $K$-free. Note also that
\[ \Orbits([H,\varphi]) = |G|/|H| \]
for a $(G,K)$-pair $(H,\varphi)$.
One easily checks that the orbit aumentations send the composition pairing
to multiplication, and consequently we obtain a $\Z$-linear functor
\[ \Orbits \colon  \Burnside \longrightarrow \Z \]
and a $\Zp$-linear functor
\[ \Orbits \colon  \CompBurnside \longrightarrow \Zp. \]
We refer to these as the \emph{orbit augmentation functors} of the
respective Burnside categories.

\subsection{The Segal conjecture}
For a finite group $G$ and a compact Lie group $K$, there is a
homomorphism of $\Z$-modules
\[
  \alpha \colon A(G,K) \longrightarrow \PtdStableMaps{BG}{BK}
\]
sending a basis element $[H,\varphi]$ to the stable map
\mbox{$\PtdClSpectrum{\varphi} \circ tr_H$}. This assignment is
functorial for $G$ and finite $K$, and so we get a functor of
$\Z$-linear categories
\[
  \alpha \colon \Burnside \longrightarrow \StableGr,
\]
which is the identity on objects.

For a finite group $G$ and a compact Lie group $K$, we have an
obvious bilinear map
\[
  A(G,1) \times A(G,K) \longrightarrow A(G,K),
\]
where $1$ denotes the trivial group, extending the map sending a
$G$-set $X$ and a $(G,K)$-bundle $Y$ to the $(G,K)$-bundle
$X \times Y$, where $G$ acts via the diagonal and $K$ acts only on
the second coordinate.

When \mbox{$K = 1$}, this gives a ring structure on $A(G) =
A(G,1)$, and the resulting ring is called the \emph{Burnside ring
of $G$}. This can also be described as the Grothendieck group of 
isomorphism classes of left $G$-sets.
For general compact Lie groups $K$, the Burnside module $A(G,K)$
becomes a module over $A(G)$ under the action described above.

The orbit augmentation functor $\Orbits$ induces a $\Z$-algebra augmentation
\[
  \epsilon \colon A(G) \longrightarrow \Z,
\]
which extends the ``counting map'' sending a $G$-set $X$ to its cardinality |X|.
We denote the augmentation ideal by $I(G)$.

The Segal conjecture states that for finite group $G$, the
homomorphism \mbox{$A(G) \to \pi^0_S(BG_+)$}, which we describe as
the composite
\[
  A(G) = A(G,1) \stackrel{\alpha}{\longrightarrow} \PtdStableMaps{BG}{B1} \cong \pi^0_S(BG_+),
\]
is a completion with respect to the ideal $I(G)$. Lewis-May-McClure
showed in \cite{LMM} that an extended version of the conjecture, describing 
the homotopy classes of stable maps between classifying spaces of finite groups,
follows from the original version. In a further generalization, May-Snaith-Zelewski 
showed in \cite{MSZ} that one can allow the target group to be a compact Lie group.
Since the Segal conjecture was
proved by Carlsson in \cite{Car}, we can state these extensions as a
theorem.
\begin{theorem}[Segal Conjecture \cite{Car,LMM,MSZ}]\label{thm:Segal}
For a finite groups $G$ and a compact Lie group $K$, the map
\[
  \alpha \colon A(G,K) \longrightarrow \PtdStableMaps{BG}{BK}
\] induces an isomorphism
\[
  \alpha_{I(G)}^{\wedge} \colon A(G,K)^{\wedge}_{I(G)} \stackrel{\cong}{\longrightarrow} \PtdStableMaps{BG}{BK},
\]
where
\[
  A(G,K)^{\wedge}_{I(G)} = \invlim{k}\left(A(G,K)/I(G)^k A(G,K)\right)
\]
 denotes the $I(G)$-completion of $A(G,K)$.
\end{theorem}

This is a magnificent result, but in general the $I(G)$-adic
completions are difficult to calculate. However, when the groups
involved are $p$-groups, the situation is simplified. For a finite
group $G$ and a compact Lie group $K$, let $\widetilde{A}(G,K)$ be
the module obtained from $A(G,K)$ by quotienting out all basis
elements of the form $[H,\triv]$, where $\triv$ is the constant
homomorphism. May-McClure noticed in \cite{MM} that when $G$ is a
$p$-group, the $I(G)$-adic topology on $\widetilde{A}(G,K)$ is
equivalent to the $p$-adic topology on $\widetilde{A}(G,K)$. 
(This is not true on $A(G,K)$.) The result is the following version of the Segal conjecture.
\begin{theorem}[Segal Conjecture \cite{Car,LMM,MM,MSZ}]\label{thm:pSegal}
Let $S$ be a finite $p$-group and $K$ be a compact Lie group. Then
the map
\[
  \alpha \colon A(S,K) \longrightarrow \PtdStableMaps{BS}{BK}
\]
induces an isomorphism
\[
  \pComp{\widetilde{\alpha}} \colon \pComp{\widetilde{A}(S,K)} \longrightarrow \StableMaps{BS}{BK}.
\]
\end{theorem}

An explanation of the target of this isomorphism is in order. Recall
that \mbox{$\PtdStable{BK} \simeq \Stable{BK} \vee
\SphereSpectrum$}. Now, the submodule of $\PtdStableMaps{BS}{BK}$
generated by the maps $\alpha([P,\triv])$, where $\triv$ is the
constant homomorphism \mbox{$P \to K$}, consists of the maps that
factor through the $\SphereSpectrum$-term of $\PtdStable{BK}$.
Therefore it is appropriate to replace $\PtdStableMaps{BS}{BK}$ with
\[
  \StableMaps{BS}{BK} \cong \PtdStableMaps{BS}{BK}/\StableMaps{BS_+}{S^0}
\]
when passing from $A(S,K)$ to $\widetilde{A}(S,K)$ to obtain the
homomorphism \mbox{$ \widetilde{\alpha} \colon \widetilde{A}(S,K)
\rightarrow \StableMaps{BS}{BK}$}. The homomorphism
$\pComp{\widetilde{\alpha}}$ is the $p$-completion of
$\widetilde{\alpha}$.

Since $\alpha$ and $\widetilde{\alpha}$ preserve composition when
restricted to finite groups, both versions of the Segal conjecture
described in this section can be formulated in terms of isomorphisms
of linear categories.

\subsection{Restrictions and transfers}\label{sec:Res&Tr}
Let $G$ be a finite group and $K$ be a compact Lie group. If $S$
is a subgroup of $G$, we can regard a $(S,K)$-pair as a
$(G,K)$-pair, or restrict a $(G,K)$-bundle to a $(S,K)$-bundle. This gives us a way to traverse between $A(S,K)$ and
$A(G,K)$, which will be useful later when we take $S$ to be a
Sylow subgroup. In this section we briefly describe these
operations and interpret them in terms of stable maps.

Given a $(S,K)$-pair $(H,\varphi)$, we can regard $H$ as a
subgroup of $G$ to get a $(G,K)$-pair $(H,\varphi)_{G}^{K}$. This
assignment extends to a homomorphism
\[
  \Phi \colon A(S,K) \longrightarrow A(G,K).
\]
It is easy to check that
\[
  [H,\varphi]_{G}^{K} = [H,\varphi]_{S}^{K} \circ [S,id]_{G}^{S}
\]
using the double coset formula. Furthermore,
 $[S,id]_{G}^{S}$ corresponds to
$[G]_{G}^{S}$, the isomorphism class of $G$, regarded as a $(S
\times G)$-set under the action \mbox{$(s,g).x = gxs^{-1}$}.

The homomorphism $\alpha$ sends $[S,id]_{G}^{S}$ to the transfer map
  \mbox{$tr_{S} \colon \PtdClSpectrum{G} \to \PtdClSpectrum{S}$}.
Thus we get an induced homomorphism
\[
  \Phi \colon \PtdStableMaps{BS}{BK} \to \PtdStableMaps{BG}{BK},~ f \mapsto f \circ tr_{S}.
\] Note that
$\alpha$ commutes with $\Phi$ by construction.

A $(G,K)$-bundle $X$, can be regarded as a $(S,K)$-bundle
$X\vert_{(S,K)}$ via restriction. This induces a
homomorphism
\[
  \Gamma \colon A(G,K) \longrightarrow  A(S,K).
\]
It is easy to see that
\[
  [X\vert_{(S,K)}] \cong [X] \circ [G]_{S}^{G},
\]
where $[X]$ denotes the isomorphism class of $X$, and
$[G]_{S}^{G}$ is the isomorphism class of $G$ regarded as a $(G\times S)$-set via the
action \mbox{$(g,s).x = gxs^{-1}$}.
 Since $[G]_{S}^{G}$
corresponds to $[S,\iota]_{S}^{G}$ in $A(S,G)$ we see
that $\Gamma$ is given by
\[
  [X] \mapsto  [X] \circ [S,\iota]_{S}^{G}.
\]
Applying $\alpha$ we get an induced homomorphism of stable maps
\[
  \Gamma \colon \PtdStableMaps{BG}{BK} \to \PtdStableMaps{BS}{BK},~ f \mapsto f \circ \PtdClSpectrum{\iota_{S}}.
\]
Note that $\alpha$ also commutes with $\Gamma$ by construction.

To accommodate the basepoint issue we need to consider the reduced
operations $\widetilde{\Phi}$ and $\widetilde{\Gamma}$ that fit into
the diagrams
\[ \xymatrix{
  {A(S,K)} \ar[0,1]^{\Phi} \ar@{->>}[d]^{\pi} & A(G,K) \ar[0,1]^{\Gamma} \ar@{->>}[d]^{\pi}  & A(S,K)  \ar@{->>}[d]^{\pi}\\
  {\widetilde{A}(S,K)} \ar[r]^{\widetilde{\Phi}} & {\widetilde{A}(G,K)} \ar[0,1]^{\widetilde{\Gamma}} & \widetilde{A}(S,K)
  }
\]
and
\[ \xymatrix{
  {\PtdStableMaps{BS}{BK}} \ar[0,1]^{\Phi} \ar@{->>}[d]^{\pi} & {\PtdStableMaps{BG}{BK}} \ar[0,1]^{\Gamma} \ar@{->>}[d]^{\pi}  & {\PtdStableMaps{BS}{BK}}  \ar@{->>}[d]^{\pi}\\
  {\StableMaps{BS}{BK}} \ar[r]^{\widetilde{\Phi}} & {\StableMaps{BG}{BK}} \ar[0,1]^{\widetilde{\Gamma}} & {\StableMaps{BS}{BK},}
  }
\]
where the vertical arrows are the obvious projection maps. The
reduced operations $\widetilde{\Phi}$ and $\widetilde{\Gamma}$
commute with $\widetilde{\alpha}$.

Since we are working $p$-locally we need to consider the
$p$-completed operations
\[
  \pComp{\widetilde{\Phi}} \colon \StableMaps{\pComp{BS}}{BK}
  \longrightarrow \StableMaps{\pComp{BG}}{BK},~ f \mapsto f \circ
  \pComp{(tr_S)}
\]
and
\[
  \pComp{\widetilde{\Gamma}} \colon \StableMaps{\pComp{BG}}{BK}
  \longrightarrow \StableMaps{\pComp{BS}}{BK},~ f \mapsto f \circ
  \pComp{(\ClSpectrum{\iota_S})}.
\]
These fit into the commutative diagram
\[ \xymatrix{
  {\widetilde{A}(S,K)} \ar[0,1]^{\widetilde\Phi} \ar[d]^{\iota_p^*\circ\widetilde{\alpha}} & \widetilde{A}(G,K) \ar[0,1]^{\widetilde\Gamma} \ar[d]^{\iota_p^*\circ\widetilde{\alpha}}  & \widetilde{A}(S,K)  \ar[d]^{\iota_p^*\circ\widetilde{\alpha}}\\
  {\StableMaps{\pComp{BS}}{BK}} \ar[r]^{\pComp{\widetilde{\Phi}}} & {\StableMaps{\pComp{BG}}{BK}} \ar[0,1]^{\pComp{\widetilde{\Gamma}}} & {\StableMaps{\pComp{BS}}{BK},}
  }
\]
where the vertical maps are the functor $\widetilde{\alpha}$
followed by precomposition with $\iota_p$, the natural inclusion of
wedge summands described in Section \ref{sec:pComps}. The naturality
of $\iota_p$ ensures the commutativity of the diagram and in fact
this is the reason why $\iota_p$ was introduced so carefully.

Notice that
\[ [G_{G}^{S}] \circ [G_{S}^{G}] \cong [G_{S}^{S}],\]
where $[G_{S}^{S}]$ denotes the isomorphism class of $G$ 
regarded as a $(S \times S)$-set via the action 
\mbox{$(s_1,s_2).x = s_1xs_2^{-1}$}. Therefore the
composite
 \[\Gamma \circ \Phi \colon A(S,K) \longrightarrow A(S,K).\]
is given by
%\[\Gamma ( \Phi (X) ) = [G_2]_{G_2}^{S_2} \circ [G_2]_{S_2}^{G_2} \circ X \circ [G_1]_{G_1}^{S_1} \circ [G_1]_{S_1}^{G_1} = [G_2]_{S_2}^{S_2} \circ X \circ [G_1]_{S_1}^{S_1}.\]
 \[\Gamma \circ \Phi (X)  = X \circ [G]_{S}^{S}.\]
This operation will be studied in Section \ref{sec:Subconjugacy}.

We also consider the opposite composite
 \[\Phi \circ \Gamma \colon A(G,K) \longrightarrow A(G,K).\]
Its corresponding operation on stable maps is given by
 \[\Phi \circ \Gamma (f) = f \circ \PtdClSpectrum{\iota_{S}} \circ tr_{S}.\]
The interesting case is when $S$ is a Sylow $p$-subgroup of $G$.
Then \mbox{$\ClSpectrum{S} \simeq \pComp{\ClSpectrum{S}}$} and the
composite \mbox{$\ClSpectrum{\iota_{S}} \circ tr_{S}$} is a homotopy
equivalence after $p$-completion. Therefore
\mbox{$\pComp{\widetilde{\Phi}} \circ \pComp{\widetilde{\Gamma}}$}
is an automorphism of $\StableMaps{\pComp{BG}}{BK}$ which factors
through $\StableMaps{BS}{BK}$. We conclude that
$\StableMaps{\pComp{BG}}{BK}$ is a direct summand of
$\StableMaps{BS}{BK}$ isomorphic to
\[
  \pComp{\widetilde{\Gamma}} \circ \pComp{\widetilde{\Phi}}\left(\StableMaps{BS}{BK}\right)
  = \pComp{\widetilde{\alpha}} \left(\pi \left(\pComp{A(S,K)} \circ [G]_S^S \right) \right).
\]
This observation gives us a strategy to prove Theorem \ref{mthm:NewSegal}. Over the next two
sections we show that $\pComp{A(S,K)} \circ [G]_S^S $ is isomorphic to $\pComp{A_p(G,K)}$ in  
a way that is compatible with $\widetilde{\alpha}_p$.

\section{Subconjugacy}\label{sec:Subconjugacy}
From now on, fix a finite group $G$ with Sylow p-subgroup $S$ and
a compact Lie group $K$. In this section we find an explicit basis
for the submodule
\[\pComp{A(S,K)} \circ [G] \subseteq \pComp{A(S,K)},\]
where \mbox{$[G]\in \pComp{A(S,S)}$} is the
isomorphism class of $G$ regarded as a \mbox{$(S\times S)$}-set under the action 
$(s_1,s_2). x = s_1xs_2^{-1}$. From
this basis it is easy to obtain a basis for
$\StableMaps{\pComp{BG}}{BK}$ and prove Theorem
\ref{mthm:NewSegal}.

The main tool we will use is a filtration induced by the following
preorder on conjugacy classes of $(S,K)$-pairs.

\begin{definition} \label{def:PairSubconjugacy}
Let $(P,\varphi)$ and $(Q,\psi)$ be two $(S,K)$-pairs. We say that
$(Q,\psi)$ is \emph{subconjugate} to $(P,\varphi)$, and write
\[ 
  (Q,\psi) \noFsub (P,\varphi),
\]
if there exist elements \mbox{$g \in G$} and \mbox{$k \in K$} such
that the following diagram commutes
\[
\begin{CD}
{Q} @ > {\psi} >> {\psi(Q)} \\
@ V c_{g} VV @ VV c_{k} V \\
{P} @> {\varphi} >> {\varphi(P).} \\
\end{CD}
\]
\end{definition}
It is clear that subconjugacy is preserved by $(S,K)$-conjugacy
and therefore we can pass to conjugacy classes.
\begin{definition} We say that a
$(S,K)$-conjugacy class of pairs $[Q,\psi]$ is subconjugate to an
$(S,K)$-conjugacy class $[P,\varphi]$ and write
\[
  [Q,\psi] \noFsub [P,\varphi]
\]
if the subconjugacy relation
\[ 
  (Q,\psi) \noFsub (P,\varphi)
\]
holds between any (and hence all) representatives of the classes.
\end{definition}

It is also clear that subconjugacy is a transitive relation and so
induces an equivalence relation as described in the following
definition.
\begin{definition}
We say that $[Q,\psi]$ is $(G,K)$-\emph{conjugate} to
$[P,\varphi]$, and write
\[ 
  [Q,\psi] \noFcon [P,\varphi],
\]
if~ \mbox{$[Q,\psi] \noFsub [P,\varphi]$}~ and~ \mbox{$[P,\varphi]
\noFsub [Q,\psi]$}. We say that $[Q,\psi]$ is \emph{strictly
subconjugate} to $[P,\varphi]$, and write
\[
  [Q,\psi] \noFstrictsub [P,\varphi],
\]
if $[Q,\psi]$ is subconjugate to $[P,\varphi]$, but not
\mbox{$(G,K)$-conjugate} to $[P,\varphi]$.
\end{definition}

Let $I$ denote the set of $(G,K)$-conjugacy classes of
$(S,K)$-pairs. This is a poset under subconjugacy and we proceed
to construct an $I$-indexed filtration of $\pComp{A(S,K)}$.
\begin{definition}
For an $(S,K)$-pair $(P,\varphi)$, let $\noFsubModule{P}{\varphi}$
denote the submodule of $\pComp{A(S,K)}$ generated by the basis
elements $[Q,\psi]$ such that
\[ 
  [Q,\psi] \noFsub (P,\varphi),
\]
and let $\noFstrictsubModule{P}{\varphi}$ denote the submodule
generated by the basis elements $[Q,\psi]$ such that
\[  
  [Q,\psi] \noFstrictsub (P,\varphi).
\]
\end{definition}

We will show that this filtration is preserved by composition with
$[G]$. First note that by the double coset formula one has
\[ 
  [G] = \sum_{x \in S\backslash G / S} [S \cap S^x, c_x ].
\]
This prompts us to make a formalization.
\begin{definition} Let $R$ be the submodule of
$\pComp{A(S,S)}$ generated by the basis elements $[P,\varphi]$,
where the homomorphism \mbox{$\varphi \colon P \to S$} is induced
by conjugation by an element in $G$.
\end{definition}
Clearly \mbox{$[G]\in R$}. Using the double coset formula, one
easily sees that $R$ is in fact a subring of $\pComp{A(S,S)}$. The
next lemma shows that the modules $\noFsubModule{P}{\varphi}$
and $\noFstrictsubModule{P}{\varphi}$ defined above are
right \mbox{$R$-modules}. They are in fact also left \mbox{$R$-modules} but this
is not needed.
\begin{lemma}\label{lem:PreservesFiltration}
The following hold for every $(S,K)$-pair $(P,\varphi)$:
\begin{itemize}
  \item[(a)] $\noFsubModule{P}{\varphi} \circ R \subseteq  \noFsubModule{P}{\varphi},$
  \item[(b)] $\noFstrictsubModule{P}{\varphi} \circ R \subseteq  \noFstrictsubModule{P}{\varphi}.$
\end{itemize}
\end{lemma}
\begin{proof}
We begin by proving part (a). It suffices to show that for every
basis element \mbox{$[Q,\psi] \noFsub [P,\varphi]$} of
$\noFsubModule{P}{\varphi}$ and every basis element $[T,c_g]$ of
$R$, one has
\[ [Q,\psi] \circ [T,c_g] \in \noFsubModule{P}{\varphi}.\]
By the double coset formula we can write
\[ [Q,\psi] \circ [T,c_g] = \sum_{x \in Q\backslash S / \lsup{T}{g}} [T\cap Q^{xg}, \psi \circ c_x \circ c_g] \]
so it suffices to note that
\[ [T\cap Q^{xg}, \psi \circ c_x \circ c_g]  \noFsub [Q,\psi] \noFsub [P,\varphi]\]
for every \mbox{$x \in S$}.

The same argument proves part (b), for if we assume that
$[Q,\psi]$ is strictly subconjugate to $[P,\varphi]$ at the
beginning, we get
\[
 [T\cap Q^{xg}, \psi \circ c_x \circ c_g] \noFsub [Q,\psi] \noFstrictsub [P,\varphi]
\]
at the end.
\end{proof}
In particular this lemma shows that the operation
\[ 
  \Gamma \circ \Phi = (-) \circ [G]  : \pComp{A(S,K)} \longrightarrow \pComp{A(S,K)}
\]
introduced in Section \ref{sec:Res&Tr} preserves the subconjugacy
filtration.

\begin{lemma} \label{lem:Ginvariant}
If the basis elements $[P,\varphi]$ and $[Q,\psi]$ of $A(S,K)$ are
$(G,K)$-conjugate, then
\[ 
  \Gamma \circ \Phi ( [P,\varphi] ) = \Gamma \circ \Phi ( [Q,\psi] ) .
\]
\end{lemma}
\begin{proof}
Clearly $(P,\varphi)$ and $(Q,\psi)$ are also $(G,K)$-conjugate
when regarded as $(G,K)$-pairs, so, using the notation of Section
\ref{sec:Res&Tr},
\[
  \Phi\left([P,\varphi]\right) = \Phi\left([Q,\psi]\right).
\]
Consequently
\[
  \Gamma \left(\Phi \left( [P,\varphi]\right) \right) = \Gamma \left(\Phi \left( [Q,\psi]\right) \right).
\]
\end{proof}

\begin{proposition} \label{prop:GAGbasis}
Pick a representative $[P_i,\varphi_i]$ for each \mbox{$i \in I$}.
The collection
\[
  C = \{ \Gamma \circ \Phi( [P_i,\varphi_i] ) \mid i \in I \}
\]
forms a $\Zp$-basis for $\Gamma \circ \Phi ( \pComp{A(S,K)}).$
\end{proposition}
\begin{proof}
By Lemma \ref{lem:Ginvariant} it is clear that $C$ spans $\Gamma
\circ \Phi ( \pComp{A(S,K)} )$, so it suffices to prove linear
independence.

For a $(S,K)$-pair $(P,\varphi)$ we have
\[
 \Gamma \circ \Phi ([P,\varphi]) = [P,\varphi] \circ [G] \in \noFsubModule{P}{\varphi}
\]
by Lemma \ref{lem:PreservesFiltration}. But note that
\[
  \Orbits{(\noFsubModule{P}{\varphi})} = \Orbits{([P,\varphi])} \cdot \Zp
\]
and
\[
 \Orbits{(\noFstrictsubModule{P}{\varphi})} = p \cdot \Orbits{([P,\varphi])} \cdot \Zp.
\]
Since
\[
 \Orbits(\Gamma \circ \Phi ([P,\varphi]) )  = \Orbits{([P,\varphi]\circ [G])} = \Orbits{([P,\varphi])} \cdot \Orbits{([G])}
\]
and $p$ does not divide \mbox{$\Orbits{([G])} = |G|/|S|$}, we deduce
that
\begin{equation} \label{eq:GPreservesLayer}
  \Gamma \circ \Phi ([P,\varphi]) \in \noFsubModule{P}{\varphi} \setminus \noFstrictsubModule{P}{\varphi}.
\end{equation}

Now, let \mbox{$c_i \in \Zp$} for each \mbox{$i \in I$} and assume
that
\begin{equation} \label{eq:LinearDependence}
  \sum_{i \in I} c_i \cdot \left( \Gamma \circ \Phi ([P_i,\varphi_i]) \right) = 0.
\end{equation}
Put
\[
  I' = \{i \in I \mid c_i \neq 0\}.
\]
If $I'$ is nonempty, then let $j$ be a maximal element of $I'$
regarded as a poset under subconjugacy. By
(\ref{eq:GPreservesLayer}) there is a $(S,K)$-pair \mbox{$(Q,\psi)
\Fcon{(G,K)} (P_j,\varphi_j)$} such that
\[
  \Coeff_{[Q,\psi]}\left( \Gamma \circ \Phi ([P_j,\varphi_j]) \right) \neq 0.
\]
On the other hand, for \mbox{$i \in I'\setminus \{j\}$}, the
maximality of $j$ implies that $[Q,\psi]$ is not subconjugate to
$(P_i,\varphi_i)$. Hence
\[ 
  \Coeff_{[Q,\psi]}\left( \noFsubModule{P_i}{\varphi_i} \right) = 0
\]
and in particular
\[ 
  \Coeff_{[Q,\psi]}\left( \Gamma \circ \Phi ([P_i,\varphi_i]) \right) = 0.
\]  
Now we get
\begin{align*}
  \Coeff_{[Q,\psi]}\left(\sum_{i \in I} c_i \cdot \left(\Gamma \circ \Phi ([P_i,\varphi_i]) \right)\right)
  &= \sum_{i \in I} c_i \cdot \Coeff_{[Q,\psi]}\left( \Gamma \circ \Phi ([P_i,\varphi_i]) \right)\\
  &= \sum_{i \in I\setminus I'} \underbrace{c_i}_{=0} \cdot \Coeff_{[Q,\psi]}\left( \Gamma \circ \Phi ([P_i,\varphi_i]) \right)\\
  &+ \sum_{i \in I'\setminus\{j\}} c_i \cdot \underbrace{\Coeff_{[Q,\psi]}\left( \Gamma \circ \Phi ([P_i,\varphi_i]) \right)}_{=0}\\
  &+ \underbrace{c_j}_{\neq 0} \cdot \underbrace{\Coeff_{[Q,\psi]}\left( \Gamma \circ \Phi ([P_i,\varphi_i]) \right)}_{\neq 0}\\
  &\neq 0,
\end{align*}
contradicting (\ref{eq:LinearDependence}). Therefore $I'$ must be
empty and we conclude that the collection is linearly independent.

\end{proof}

\section{$p$-completed Segal conjectures} \label{sec:NewSegal}
In this section we state and prove the new variants of the Segal
conjecture promised in the introduction. We begin by describing the
group of homotopy classes of stable maps from the $p$-completed 
classifying space of a finite group to the classifying space of a compact
Lie group. The map $\pComp{(\iota_p^* \circ \widetilde{\alpha}_p)}$ in
the statement of the theorem is the algebraic $p$-completion of the
composite
\[ \widetilde{A}(G,K) \xrightarrow{\widetilde{\alpha}_p} \StableMaps{BG}{BK} \xrightarrow{- \circ \iota_p} \StableMaps{\pComp{BG}}{BK},\]
where \mbox{$\iota_p \colon \pComp{\ClSpectrum{G}} \hookrightarrow
\ClSpectrum{G}$} is the natural wedge-summand inclusion discussed in
Section \ref{sec:pComps}. Note that the target
$\StableMaps{\pComp{BG}}{BK}$ is $p$-complete because it is a direct
summand of $\StableMaps{BS}{BK}$ which is $p$-complete by the Segal
Conjecture.

\begin{theorem} \label{thm:NewSegal}
For a finite group $G$ and a compact Lie group $K$, the map
\[
  \pComp{(\iota_p^* \circ \widetilde{\alpha}_p)} \colon \pComp{\widetilde{A}_p(G,K)} \longrightarrow\StableMaps{\pComp{BG}}{BK}
\]
is an isomorphism of
$\Zp$-modules, natural in $G$ and finite $K$.
\end{theorem}
\begin{proof}
Let $S$ be a Sylow subgroup of $G$, and let $I$ be the set of 
$(G,K)$-conjugacy classes of $(S,K)$-pairs. For each \mbox{$i \in
I$}, pick a representative $(P_i,\varphi_i)$ and consider the
submodule \mbox{$M_I \subseteq \pComp{A(S,K)}$} generated by the
collection $\{[P_i,\varphi_i] \mid i \in I \}$. Since every
$p$-subgroup of $G$ is conjugate to a subgroup of $S$, the
homomorphism \mbox{$\Phi \colon A(S,K) \to A(G,K)$}
described in Section \ref{sec:Res&Tr} induces an isomorphism
\[
  \pComp{\Phi}  : M_I \stackrel{\cong}{\longrightarrow} \pComp{A_p(G,K)}.
\]
Letting $\widetilde{M}_I$ denote the image of $M_I$ in
$\pComp{\widetilde{A}(S,K)}$, this descends to an isomorphism
\[
  \pComp{\widetilde{\Phi}}  : \widetilde{M}_I \stackrel{\cong}{\longrightarrow} \pComp{\widetilde{A}_p(G,K)}.
\]
Proposition \ref{prop:GAGbasis} implies that there is an
isomorphism
\[
  \pComp{\widetilde{\Gamma}} \circ \pComp{\widetilde{\Phi}} \colon \widetilde{M}_I \stackrel{\cong}{\longrightarrow} \pComp{\widetilde{\Gamma}}\circ\pComp{\widetilde{\Phi}}(\widetilde{M}_I)
\]
and that
\[
  \pComp{\widetilde{\Gamma}} \circ \pComp{\widetilde{\Phi}}(\widetilde{M}_I) =   \pComp{\widetilde{\Gamma}} \circ \pComp{\widetilde{\Phi}} (\pComp{\widetilde{A}(S,K)}).
\]
By the Segal conjecture, and the commutativity of $\alpha$ with
$\pComp{\widetilde{\Gamma}}$ and $\pComp{\widetilde{\Phi}}$, there
is an isomorphism
\[
 \pComp{\widetilde{\alpha}} \colon \pComp{\widetilde{\Gamma}} \circ \pComp{\widetilde{\Phi}} (\pComp{\widetilde{A}(S,K)}) \stackrel{\cong}{\longrightarrow}  \pComp{\widetilde{\Gamma}} \circ \widetilde{\Phi} \left(\StableMaps{BS}{BK}\right)
\]
Since the composition $\ClSpectrum{\iota_{S}} \circ tr_{S}$ is a
homotopy self-equivalence of $\pComp{\ClSpectrum{G}}$, the
homomorphisms $\pComp{\widetilde\Phi}$ and $
\pComp{\widetilde{\Gamma}}$ induce an automorphism
\[
 \pComp{\widetilde\Phi} \circ  \pComp{\widetilde{\Gamma}} \colon \StableMaps{\pComp{BG}}{BK} \stackrel{\cong}{\longrightarrow}
 \StableMaps{\pComp{BG}}{BK}.
\]
We conclude that $\pComp{\widetilde{\Phi}}$ is a surjection, so
\[
 \StableMaps{\pComp{BG}}{BK} = \pComp{\widetilde{\Phi}} \left( \StableMaps{BS}{BK}
 \right),
\]
and that $\pComp{\widetilde{\Gamma}}$ is an injection. Hence there
is an isomorphism
\[
 \pComp{\widetilde{\Gamma}} \colon \StableMaps{\pComp{BG}}{BK} \stackrel{\cong}{\longrightarrow}  \pComp{\widetilde{\Gamma}} \circ \pComp{\widetilde{\Phi}} \left(\StableMaps{BS}{BK}\right).
\]

Now we have a commutative diagram
\[
\xymatrix{
  {\widetilde{M}_I} \ar[0,1]^{\Gamma\circ \Phi \phantom{\Phi\Phi\Phi}}_{\cong\phantom{\Phi\Phi\Phi}} \ar[1,1]^{\Phi}_{\cong}  & {\Gamma\circ\Phi(\widetilde{M}_I)} \ar[0,1]^{= \phantom{\Phi\Phi\Phi}} & {\Gamma\circ\Phi(\pComp{\widetilde{A}(S,K)})} \ar[0,1]^{\pComp{\widetilde\alpha}\phantom{\Phi}}_{\cong\phantom{\Phi}} & {\Gamma\circ \Phi(\StableMaps{BS}{BK})} \\
  & {\pComp{\widetilde{A}_p(G,K)}}\ar[0,2]^{\pComp{(\iota_p^* \circ \widetilde{\alpha}_p)}}\ar[-1,0]^{\Gamma} & &  {\StableMaps{\pComp{BG}}{BK},} \ar[-1,0]^{\Gamma}_{\cong}
  }
\]
where we have written $\Phi$ and $\Gamma$ instead of
$\pComp{\widetilde{\Phi}}$ and $\pComp{\widetilde{\Gamma}}$.
Concentrating first on the left side of this diagram, we see that
since $\Phi$ and \mbox{$\Gamma \circ \Phi$} are both isomorphisms,
$\Gamma$ must be an isomorphism. Turning our attention to the right
side of the diagram, we see that the bottom arrow forms the bottom
part of a rectangular commutative diagram where all the other arrows
are isomorphism. Hence the bottom arrow is itself an isomorphism.
\end{proof}

Using Lemma \ref{lem:pCompFactors} we obtain an immediate reformulation.
%We can reformulate this theorem as follows.
\begin{corollary} \label{cor:NewSegalComp}
For a finite group $G$ and a compact Lie group $K$, the map
\[
  \pComp{(\widetilde{\alpha}_{p})} \colon \pComp{\widetilde{A}_p(G,K)} \longrightarrow \pComp{\StableMaps{BG}{BK}} \simeq \CompStableMaps{BG}{BK},
\]
is an isomorphism of $\Zp$-modules, natural in $G$ and finite $K$.
\end{corollary}

The naturality property in Corollary \ref{cor:NewSegalComp} can also be stated as follows.
\begin{corollary} \label{Cor:NewSegalCat}
The functor
\[
  \pComp{(\widetilde{\alpha}_p)} \colon \CompBurnsidep \to \CompStableGr
\]
is an isomorphism of $\Zp$-linear categories.
\end{corollary}

Finally we can collect our results for different primes $p$ and
prove Theorem \ref{mthm:NewSegalSum} of the introduction.
\begin{theorem} \label{thm:NewSegalSum}
For a finite group $G$ and a compact Lie group $K$, the map
\[
  \bigoplus_q \qComp{(\iota_q^*\circ \widetilde{\alpha}_q)} \colon \bigoplus_q \qComp{\widetilde{A}_q(G,K)} \to \bigoplus_q \StableMaps{\qComp{BG}}{BK} \cong \StableMaps{BG}{BK},
\]
where the sums run over all primes $q$, is an isomorphism of
$\Z$-modules, natural in $G$ and finite $K$.
\end{theorem}
\begin{proof}
The stable splitting \mbox{$\ClSpectrum{G} \simeq \bigvee_{q}
\qComp{\ClSpectrum{G}}$} induces a splitting
\[
  \StableMaps{BG}{BK} \cong \bigoplus_{q} \StableMaps{\qComp{BG}}{BK}.
\]
The result now follows from Theorem \ref{thm:NewSegal}.
\end{proof}
One can also restrict this theorem to finite groups and formulate a
corollary about an isomorphism of $\Z$-linear categories, but we
refrain from doing so.

\section{Alternative formulations} \label{sec:Alternative}
The basepoint issue mentioned before in this paper is more a
technical nuisance for formulating statements than an actual problem in proving them. 
The main issue with the added basepoint is that for a finite group $G$, the
suspension spectrum of the $p$-completion of $BG_+$ is
$\Stable{(\pComp{(BG_+)})} \simeq \pComp{\ClSpectrum{G}}\vee
\SphereSpectrum$, which is not $p$-complete as $\SphereSpectrum$ is
not $p$-complete. Consequently, for a compact Lie group $K$, the
module $\StableMaps{\pComp{(BG_+)}}{BK_+}$ need not be $p$-complete
as maps factoring through the sphere spectrum can contribute a
non-completed part.

In this paper we have so far opted to remove the offending sphere
spectrum, at the cost of disregarding maps factoring
through the sphere spectrum, which results in introducing
the aesthetically unpleasant quotients of double Burnside modules.
However, there are other ways to smooth over the basepoint issue. In
this section we consider two alternative approaches. The first is to
$p$-complete the added sphere spectrum rather than to remove it, and
thus force everything in sight to be $p$-complete. The second
approach is to deal separately with the failure of
$\StableMaps{\pComp{(BG_+)}}{BK_+}$ to be $p$-complete, which turns
out to be precisely one copy of $\Z$ corresponding to maps from one
sphere spectrum to the other.

For each approach we give the corresponding versions of of Theorem
\ref{thm:Segal} and state the generalizations to $p$-completed
classifying spaces of finite groups. These generalizations are
obtained using the same filtration argument already presented and so
the proofs are omitted.

\subsection{$p$-completed suspension spectra}
Instead of looking at the suspension spectrum of the $p$-completion
of $BG_+$, for a finite group $G$, we can look at the $p$-completion
of the suspension spectrum of $BG_+$. Thus we obtain the
$p$-complete spectrum \mbox{$\pComp{(\PtdClSpectrum{G})} \simeq
\pComp{\ClSpectrum{G}} \vee \pComp{\SphereSpectrum}$}. Taking this
approach for a finite $p$-group $S$ and a compact Lie group $K$, the
appropriate version of Theorem \ref{thm:Segal} is that the map
\[
  \pComp{\alpha} \colon \pComp{A(S,K)}
  \stackrel{\cong}{\longrightarrow} [\pComp{(\PtdClSpectrum{S})},\pComp{(\PtdClSpectrum{K})}]
\]
is an isomorphism of $\Zp$-modules. Using the same transfer and
filtration arguments as in the proof of Theorem \ref{thm:NewSegal}
we obtain the following generalization.
\begin{theorem} For a finite group $G$ and a compact Lie group $K$,
the map
\[
  \pComp{(\alpha_p)} \colon \pComp{A_p(G,K)}
  \stackrel{\cong}{\longrightarrow} [\pComp{(\PtdClSpectrum{G})},\pComp{(\PtdClSpectrum{K})}]
\]
is an isomorphism of $\Zp$-modules
\end{theorem}

While this approach has the advantage of retaining more information
from the double Burnside module $A_p(G,K)$, it has the drawback that
 we obtain
no analogue of Theorem \ref{thm:NewSegalSum} as $\PtdClSpectrum{G}$
is not a wedge sum of its completions (since the sphere spectrum is
not a torsion spectrum).

\subsection{Precise approach}
The statement of the Segal Conjecture for stable maps from
classifying spaces of $p$-groups presented in Theorem
\ref{thm:Segal} is actually a simplified version of a more detailed
result obtained by May-McClure in \cite{MM}.

For a finite group $G$ with Sylow subgroup $S$ and a compact Lie
group $K$, let $I_p(G,K)$ be the kernel of the map $\Orbits \colon
A_p(G,K) \to \Z$. As a $\Z$-module, $I_p(G,K)$ has basis $\{
[P,\varphi] - |S|/|P| \cdot [S,\triv] \}$. There is a splitting
\mbox{$A_p(G,K) = \Z \oplus I_p(G,K),$} where the $\Z$ term is
generated by $[S,\triv]$. Since \mbox{$ I(G) \cdot A(G,K) = I(G,K)
$} this induces a splitting \mbox{$A_p(G,K)^{\wedge}_{I(G)} = \Z
\oplus I_p(G,K)^{\wedge}_{I(G)}$}.

May-McClure showed that for a finite $p$-group $S$, the $I(S)$-adic
topology on $I(S,K)$ coincides with the $p$-adic topology. Therefore
the Segal Conjecture says that the map
\[
  (\alpha_p)^{\wedge}_{I(S)} \colon \Z \oplus \pComp{I(S,K)} \longrightarrow \PtdStableMaps{BS}{BK}
\]
is an isomorphism. Again, using the transfer and filtration
arguments from Sections \ref{sec:Subconjugacy} and
\ref{sec:NewSegal}, one can prove the following generalization.
\begin{theorem}
For a finite group $G$ and a compact Lie group $K$, the map
\[
  (\alpha_p)^{\wedge}_{I(G)} \colon \Z \oplus \pComp{I_p(G,K)} \longrightarrow \StableMaps{\pComp{(BG_+)}}{BK_+}
\]
is an isomorphism of $\Z$-modules.
\end{theorem}

Collecting this result for different primes one can obtain an isomorphism
of $\Z$-modules.
\[
  \Z \oplus \bigoplus_p \pComp{I_p(G,K)} \stackrel{\cong}{\longrightarrow} \PtdStableMaps{BG}{BK}.
\] 
Unlike the isomorphism in \ref{thm:NewSegalSum} this isomorphism is not natural in $G$.

Taking $K = 1$ one recovers a special case of Minami's description of the $I(G)$-adic completion of 
the Burnside ring in \cite{Min} (he allows $G$ to be compact). It should be noted that even though one
can regard $A(G,K)$ as a submodule of $A(K\times G)$ when $K$ is finite, the results in this paper do not
follow from Minami's results, as his result involves the $I(K\times G)$-adic completion whereas we are 
interested in the $I(G)$-adic completion, and the two actions are hard to reconcile.

\section{Decompositions} \label{sec:Decomposition}
In this section we address the following question. Given a
$(G,K)$-pair $(H,\varphi)$, what is the element in
$\pComp{\widetilde{A}_p(G,K)}$ corresponding to the $p$-completed stable
map $\widetilde{\alpha}([H,\varphi])$? In other words, we describe the homomorphism
\[\widetilde{\pi}_p \colon \widetilde{A}(G,K) \xrightarrow{\pComp{(-)}\circ\widetilde\alpha} \CompStableMaps{BG}{BK} \xrightarrow{(\pComp{(\widetilde{\alpha}_p)})^{-1}} \pComp{\widetilde{A}_p(G,K)}.\]

For technical reasons it is more convenient to first describe the homomorphism 
\[{\pi}_p \colon {A}(G,K) \xrightarrow{\pComp{(-)}\circ \alpha} [\pComp{(\PtdClSpectrum{G})},\pComp{(\PtdClSpectrum{K})}] \xrightarrow{(\pComp{(\alpha_p)})^{-1}} \pComp{A_p(G,K)}.\]
and then interpret the results for the homomorphism $\widetilde{\pi}$ obtained 
by removing the additional sphere spectrum.

Let $1_p \in  \pComp{A_p(G,G)}$ be the pre-image of the identity of $\pComp{(\PtdClSpectrum{G})}$ under the isomorphism
 \[ \pComp{(\alpha_p)} \colon \pComp{A_p(G,G)} \stackrel{\cong}{\longrightarrow} [\pComp{(\PtdClSpectrum{G})},\pComp{(\PtdClSpectrum{G})}]. \]
This immediately gives us a way to describe $\pi_p$, for if \mbox{$X \in A(G,K)$} then we can regard $X \circ 1_p$ as an element of $\pComp{A_p(G,K)}$, and we have
 \[\pComp{(\alpha_p)}(X \circ 1_p) = \pComp{\alpha(X)} \circ \pComp{(\alpha_p)}(1_p) = \pComp{\alpha(X)} \circ id = \pComp{\alpha(X)}.\]
Therefore
 \[\pi_p (X) = X \circ 1_p.\]
We proceed to determine $1_p$ using
filtration methods similar to those in Section
\ref{sec:Subconjugacy}. Since many details are similar, they are
left to the reader.

Fix a Sylow
subgroup $S$ of $G$. Like in Section \ref{sec:Subconjugacy}, we
say that a $(G,G)$-pair $(F,\psi)$ is subconjugate to a
$(G,G)$-pair $(H,\varphi)$, and write \mbox{$(F,\psi) \noFsub
(H,\varphi) $} if there exist elements \mbox{$g_1,g_2 \in G$}
making the following diagram commute
\[
\begin{CD}
{F} @ > {\psi} >> {\psi(F)} \\
@ V c_{g_1} VV @ VV c_{g_2} V \\
{H} @> {\varphi} >> {\varphi(H).} \\
\end{CD}
\]
This time the induced equivalence relation is just
$(G,G)$-conjugacy, and so subconjugacy induces an order relation
on $(G,G)$-conjugacy classes of $(G,G)$-pairs.

For a basis element $[P,\varphi]$ of $A_p(G,G)$, let
$\noFsubModule{P}{\varphi}$ denote the submodule of $A_p(G,G)$
generated by basis elements that are subconjugate to
$[P,\varphi]$, and similarly let $\noFstrictsubModule{P}{\varphi}$
denote the submodule generated by strictly subconjugate basis
elements. We make the following key observation, which is proved
just like Proposition
\ref{prop:GAGbasis}.
\begin{lemma} \label{lem:SPreservesFiltration}
For every $(G,G)$-pair $(P,\varphi)$, where $P$ is a $p$-group,
 \[[P,\varphi] \circ [S,\iota_S] \in \noFsubModule{P}{\varphi} \setminus \noFstrictsubModule{P}{\varphi}.\]
\end{lemma}
%\begin{proof}
%One easily proves that \mbox{$[P,\varphi] \circ [S,\iota_S] \in
%\noFsubModule{P}{\varphi}$} using the double coset formula. To
%prove that \mbox{$[P,\varphi] \circ [S,\iota_S] \notin
%\noFstrictsubModule{P}{\varphi}$} one notes that
 %\[\Orbits \left(\noFstrictsubModule{P}{\varphi} \right) = p \cdot \frac{|G|}{|P|} \cdot \Zp,\]
%while
% \[\Orbits ([P,\varphi] \circ [S,\iota_S]) = \frac{|G|}{|P|} \cdot \frac{|G|}{|S|},\]
%and $p$ does not divide $\frac{|G|}{|S|}$.
%\end{proof}
Let $R$ denote the subring of $A_p(G,G)$ generated by basis
elements of the form $[P,\iota_P]$, where $\iota_P$ denotes the
inclusion \mbox{$P \leq G$}. Note that \mbox{$R =
\noFsubModule{S}{\iota}$}. The above lemma allows us to prove the
following.
\begin{proposition}
$1_p \in R.$
\end{proposition}
\begin{proof}
Take a maximal $[Q,\psi]$ such that \mbox{$\Coeff_{[Q,\psi]}(1_p)
\neq 0$}. By Lemma \ref{lem:SPreservesFiltration} we have
 \[\Coeff_{[Q,\psi]} ([Q,\psi] \circ [S,\iota_S]) \neq 0\]
(otherwise \mbox{$[Q,\psi] \circ [S,\iota_S] \in
\noFstrictsubModule{Q}{\psi}$}).

If $[P,\varphi]$ is a basis element different from $[Q,\psi]$ such
that \mbox{$\Coeff_{[P,\varphi]} (1_p) \neq 0$}, then, by
maximality, $[Q,\psi]$ is not subconjugate to $[P,\varphi]$ and
hence \mbox{$\Coeff_{[Q,\psi]}(X) = 0$} for all \mbox{$X
\in \noFsubModule{P}{\varphi}$}. In particular,
 \[\Coeff_{[Q,\psi]}([P,\varphi] \circ [S,\iota_S]) = 0.\]

Now we deduce that
\begin{align*}
  \Coeff_{[Q,\psi]}\left(1_p \circ [S,\iota_S]\right)
  &= \Coeff_{[Q,\psi]}\left(\sum_{[P,\varphi]} \Coeff_{[P,\varphi]}(1_p) \cdot [P,\varphi] \circ [S,\iota_S]\right) \\
  &= \sum_{[P,\varphi]} \Coeff_{[P,\varphi]}(1_p) \cdot \Coeff_{[Q,\psi]}\left([P,\varphi] \circ [S,\iota_S]\right)\\
  &= \Coeff_{[Q,\psi]}(1_p) \cdot \Coeff_{[Q,\psi]}\left([Q,\psi] \circ  [S,\iota_S]\right)\\
  &\neq 0.
\end{align*}
On the other hand, we have \mbox{$1_p \circ [S,\iota_S] =
[S,\iota_S]$}. Therefore we must have \mbox{$[Q,\psi]=
[S,\iota_S]$}. We deduce that $[S,\iota_S]$ is the unique maximal element such that \mbox{$\Coeff_{[S,\iota_S]}(1_p) \neq 0$}, and hence \mbox{$1_p \in
\noFsubModule{S}{\iota_S} = R$}.
\end{proof}

Let $n$ be the number of conjugacy classes of $p$-subgroups
of $G$ and pick one representative $P_i$ for each conjugacy class,
labelled from $1$ to $n$ so that \mbox{$P_1 = S$}, and $i \geq j$
if $P_j$ is conjugate to a subgroup of $P_i$. Since \mbox{$1_p \in
R$}, we can write
 \[1_p = \sum_{j = 0}^n a_j \cdot [P_j,\iota_{P_j}].\]

For subgroups $H$ and $F$ of $G$, let $N_G(F,H)$ denote the
\emph{transporter}
 \[N_G(K,H) = \{ g \in G \mid K^g \leq H \}.\]

\begin{proposition} \label{prop:TriangularMatrix}
The multiplicative identity of $\pComp{A_p(G,G)}$ is given by
 \[ 1_p = \sum_{j = 1}^n a_j \cdot [P_j,\iota_{P_j}],  \]
where the coefficients $a_j$ satisfy the equations
\begin{equation} \label{eq:TriangularMatrix}
\sum_{j=1}^n a_j \cdot \frac{|N_G(P_i,P_j)|}{|P_j|} = 1
\end{equation}
for \mbox{$i = 1,\ldots,n$}.
\end{proposition}
\begin{proof}
In this proof we will denote all inclusions of groups by $\iota$.
For every $i$ we have \mbox{$1_p \circ [P_i,\iota] =
[P_i,\iota]$}. In particular, \mbox{$\Coeff_{[P_i,\iota]}(1_p \circ
[P_i,\iota]) = 1$}. By the double coset formula, we have
 \[[P_j,\iota] \circ [P_i,\iota] = \sum_{x \in P_j \backslash G / P_i} [P_i \cap P_j^x,c_x].\]
It is easy to check that
  \mbox{$[P_i\cap P_j^x,c_x] = [P_i \cap P_j^x,\iota] $},
and
 \[\Coeff_{[P_i,\iota]}([P_i \cap P_j^x,\iota])  = \begin{cases}
                                                       {1} & \text{if $\lsup{P_i}{x} \leq  P_j$},\\
                                                       {0} & \text{otherwise}.
                                                       \end{cases}
 \]
Noting that \mbox{$ \lsup{P_i}{x} \leq  P_j$} if and only if \mbox{$x
\in N_G(P_i,P_j)$}, we see that
\[
  \Coeff_{[P_i,\iota]}([P_j,\iota] \circ [P_i,\iota]) = |P_j \backslash N_G(P_i,P_j) / P_i|
\]
If \mbox{$x \in N_G(P_i,P_j)$} and \mbox{$g \in P_i$}, then
there is a \mbox{$h \in P_j$} such that \mbox{$xg = hx$}.
Therefore
 \[P_j \backslash N_G(P_i,P_j)/ P_i \cong P_j \backslash N_G(P_i,P_j),\]
and
  \[\Coeff_{[P_i,\iota]}([P_j,\iota] \circ [P_i,\iota]) = \frac{|N_G(P_i,P_j)|}{|P_j|}.\]
Now we get
\begin{align*}
    1 &= \Coeff_{[P_i,\iota]}(1_p \circ [P_i,\iota])\\
      &= \Coeff_{[P_i,\iota]}(\sum_{j = 0}^n a_j \cdot [P_j,\iota] \circ [P_i,\iota])\\
      &= \sum_{j = 1}^n a_j \cdot \Coeff_{[P_i,\iota]}([P_j,\iota] \circ [P_i,\iota])\\
      &= \sum_{j = 1}^n a_j \cdot \frac{|N_G(P_i,P_j)|}{|P_j|}
\end{align*}

\end{proof}

\begin{remark}
Since \mbox{$|N_G(P_i,P_j)| = 0$} when \mbox{$i < j$}, the
equations in (\ref{eq:TriangularMatrix}) constitute a lower
triangular matrix equation. The elements on the diagonal are
\mbox{$|N_G(P_i)/P_i| \neq 0$}, so the matrix is of maximal rank
and the equations suffice to determine $1_p$ uniquely. We do not
need to prove that the equation is solvable (this could be an
issue because the elements on the diagonal need not be units in
$\Zp$) since we already know that $\pComp{A_p(G,G)}$ has a
multiplicative identity.
\end{remark}

We now proceed to describe $\pi_p$ in terms of fixed-point sets.
For every $(G,K)$-pair $(H,\psi)$, there is a homomorphism
\[ \EulerSub{[H,\psi]} \colon A(G,K) \longrightarrow \Z , \]
depending only on the conjugacy class of $(H,\psi)$, which sends a $(G,K)$-bundle $X$ to
 \[ \EulerSub{[H,\psi]} (X) = \Euler\left(W(\Graph{H}{\psi})\backslash X^{\Graph{H}{\psi}}\right), \]
where $\Euler$ is the Euler characteristic, $X^{\Graph{H}{\psi}}$ is the fixed-point space, and
 \[ W(\Graph{H}{\psi}) = N_{K\times G}(\Graph{H}{\psi})/\Graph{H}{\psi}. \]
Note that for a $(G,K)$-pair $(H',\psi')$ we have
\[ W(\Graph{H}{\psi})\backslash((K\times G)/\Graph{H'}{\psi'})^{\Graph{H}{\psi}} \cong N_{K\times G}(\Graph{H}{\psi})\backslash N_{K\times G}(\Graph{H}{\psi},\Graph{H'}{\psi'})/\Graph{H'}{\psi'} \]
and that $N_{K\times G}(\Graph{H}{\psi})\backslash N_{K\times G}(\Graph{H}{\psi},\Graph{H'}{\psi'})$ is finite, so
we in fact have
\[ \EulerSub{[H,\psi]} (X) = \left|W(\Graph{H}{\psi})\backslash X^{\Graph{H}{\psi}}\right| \]
for a finite principal $(G,K)$-bundle $X$ over a finite $G$-set. 

Let $C_p(G,K)$ be the set of conjugacy classes of $(G,K)$-pairs $(P,\psi)$ where $P$ is a $p$-group.   
\begin{proposition} \cite{Dieck}
The homomorphism
 \[ \EulerSub{p} \colon A_p(G,K) \longrightarrow \prod_{ C_p(G,K)} \Z,~~~ 
 %  \hphantom{\EulerSub{[H,\psi]}}\]
 %\[ \hphantom{ \EulerSub{C_p(G,K)} \colon A_p(G,K)  } 
 X \longmapsto \prod_{C_p(G,K)} \EulerSub{[P,\psi]}(X)  \]
is an injection.
\end{proposition}

We will describe $\pi_p$ in terms of this embedding. To this end we need the following reformulation of a result of Benson--Feschbach.
\begin{lemma}\cite{BF}
For $X \in A(G,K)$, a $(G,K)$-pair $(H,\psi)$ and a subgroup $F \leq G$, we have
\[ \EulerSub{[H,\psi]}(X \circ [F,\iota_F]_G^G) = \frac{|N_{G}(H,F)|}{|F|} \EulerSub{[H,\psi]}(X)  .\]
In particular $\EulerSub{[H,\psi]}(X \circ [F,\iota_F]_G^G) = 0$ if $H$ is not subconjugate to $F$.
\end{lemma}
\begin{proof} This follows from Proposition 3.1 in \cite{BF} (which is actually proved in the case where $K$ is finite and without taking the quotient by $W(\Graph{H}{\psi})$ but the same argument works in this setting).
\end{proof}

\begin{theorem}
The homomorphism $\pi_p$
%\[{\pi}_p \colon : A(G,K) \longrightarrow \pComp{A_p(G,K)}\]
sends $X \in A(G,K)$ to the unique element $X_p \in \pComp{A_p(G,K)}$ such that
\[ \Euler_p(X_p) = \Euler_p(X). \]
In other words, $\pi_p$ is the unique homomorphism fitting into the commutative diagram.
\[ \xymatrix{
  A(G,K) \ar[rr]^{\pi_p} \ar[dr]^{\chi_p} && \pComp{A_p(G,K)} \ar@{>->}[dl]_{\chi_p}\\
  & \prod\limits_{C_p(G,K)}\Zp. &   
 }
\]
The homomorphism $\widetilde{\pi}_p$ is the homomorphism obtained from $\pi_p$ by quotienting out all trivial basis elements.
\end{theorem}
\begin{proof} 
For $X \in A(G,K)$ and a $(G,K)$-pair $(Q,\psi)$ with $Q$ a $p$-group we have 
\begin{align*}
  \EulerSub{[Q,\psi]} (\pi_p(X)) 
  &= \EulerSub{[Q,\psi]}(X\circ 1_p)\\
  &= \sum_{j = 0}^n a_j \EulerSub{[Q,\psi]}( X \circ [P_j,\iota_{P_j}] ) \\
  &= \sum_{j = 0}^n a_j \frac{|N_G(Q,P_j)|}{|P_j|}  \EulerSub{[Q,\psi]}( X ) \\
  &= \EulerSub{[Q,\psi]}( X ),
\end{align*}
where the last step follows from Proposition \ref{prop:TriangularMatrix} since $Q$ is conjugate to one of the groups $P_i$.
\end{proof}

Letting $p$ run over all primes we obtain the following corollary, from which Theorem \ref{mthm:Kernel} in the introduction follows.
\begin{corollary} \label{cor:Kernel} The kernel of the map 
 \[\alpha \colon A(G,K) \longrightarrow \StableMaps{BG_+}{BK_+} \cong A(G,K)^{\wedge}_{I(G)} \] 
is the submodule consisting of elements $X$ such that
$ \EulerSub{[H,\psi]}(X) = 0 $ for all $(G,K)$-pairs $[H,\psi]$ where the order of $H$ is a prime power.
\end{corollary}

\begin{remark} \label{rem:DropW} If $K$ is finite one may replace the morphism $\EulerSub{[H,\psi]}$ with the morphism sending $X \mapsto |X^{\Graph{H}{\psi}}|$ throughout this section and the results still hold true. 
\end{remark}

\end{document}